\documentclass[12pt]{amsart}
\usepackage{amssymb}
\usepackage{amscd}

\theoremstyle{plain}
\newtheorem{thm}{Theorem}

\newtheorem{prop}[thm]{Proposition}

\newtheorem{Def}{Definition}

\theoremstyle{definition}

\newcommand{\Hom}{\operatorname{Hom}}
\newcommand{\Aut}{\operatorname{Aut}}
\newcommand{\Inn}{\operatorname{Inn}}

\newcommand{\Tor}{\operatorname{Tor}}
\newcommand{\Isom}{\operatorname{Isom}}
\newcommand{\Id}{\operatorname{Id}}
\newcommand{\RP}{{\mathbb R}{\mathbb P}}
\newcommand{\R}{{\mathbb R}}
\newcommand{\cH}{{\mathcal H}}
\newcommand{\cC}{{\mathcal C}}
\newcommand{\cG}{{\mathcal G}}

\newcommand{\End}{\operatorname{End}}
\newcommand{\M}{\operatorname{M}}

\newcommand{\Br}{\operatorname{Br}}
\newcommand{\C}{{\mathbb C}}
\newcommand{\HH}{{\mathbb H}}
\newcommand{\Rbox}{{\begin{flushright} $\blacksquare$ \end{flushright}}}

\author{Finlay Thompson}
\thanks{Victoria University of Wellington.}
\address{\hskip-\parindent 
Finlay Thompson\\
VUW \\
Wellington, New Zealand}
\email{Finlay.Thompson@mcs.vuw.ac.nz}
\date{\today}

\begin{document}
\title[Quaternionic Gerbes]{{\large Introducing Quaternionic Gerbes.}} 
\begin{abstract} The notion of a quaternionic gerbe is
  presented as a new way of bundling algebraic structures over a four
  manifold. The structure groupoid of this fibration is described in
  some detail. The Euclidean conformal group $\R^+SO(4)$ appears
  naturally as a (non-commutative) monoidal structure on this
  groupoid. Using this monoidal structure we indicate the 
  existence of a canonical quaternionic gerbe associated to a conformal
  structure on a four manifold. \end{abstract} 
\maketitle

It is natural to think that quaternionic algebra and four dimensional
geometry should be closely linked. Certainly complex algebra and
analysis provide indispensable tools for exploring two dimensional
Riemaniann geometry.

However, despite many attempts, quaternionic algebra has not been
usefully applied to the differential geometry of four
manifolds.\footnote{Except perhaps Atiyah's notes on solutions to the
Yang-Mills equations on the four sphere, \cite{Ati}} The most commonly held
view is that quaternionic algebra is too rigid to be useful
in studying four manifolds. It is generally assumed that the natural
setting for quaternionic differential geometry is
hyperK\"{a}hler or hypercomplex. \cite{Joy}  

The purpose of this talk/article is to present the notion of a
\textbf{quaternionic gerbe}, and to demonstrate that they appear
naturally as a quaternionic algebraic structure on four
manifolds. This work appears as part of an effort to realize the goal
of ``doing four dimensional geometry and topology with quaternionic
algebra.''

Although quaternionic structures are defined \cite{Sal} for all $4n$
dimensional manifolds, the basic structures and difficulties are
already present in only 
four dimensions. The notion of ``quaternionic curve'' has been
equated with that of a ``self dual conformal'' structure.\cite{Besse}
Note that even this class of manifolds is strictly larger than the
hyperK\"{a}hler manifolds. Here we restrict our attention to smooth
oriented four manifolds, including hyperK\"{a}hler and self dual
conformal manifolds. 

It is proposed that a ``quaternionic structure'' on a four manifold
is essentially a Euclidean conformal structure. This compares favourably
with the two dimensional case where fixing a complex structure is
equivalent to fixing a conformal structure.

\section{The Problem.}

The most obvious definition of a quaternionic structure on a four
manifold $M$ requires the existence of three integrable complex
structures, $I,J,K\in\End(TM)$, such that,
$$ I^2 = J^2 =K^2 = IJK = -1 .$$

In terms of holonomy, this implies a reduction of the
frame bundle's structure group to $\HH^*$, the group of unit
quaternions. Note that $\HH^* = GL(1,\HH)$, which generalises the complex
case in an obvious way. 

The problem comes when we consider Berger's list \cite{Besse} of
holonomy groups for Riemannian manifolds:
\begin{eqnarray*}
\mbox{real} &&O(n),\quad SO(n),\\
\mbox{complex}&& U(n),\quad SU(n), \\
\mbox{quaternionic}&& Sp(n)\cdot Sp(1),\quad Sp(n) \\
\mbox{exceptional}&& G_2,\quad Spin(7)
\end{eqnarray*}

The \textbf{quaternionic-K\"{a}hler} series $Sp(n)\cdot Sp(1)$ is clearly
related to quaternionic algebraic structures, however it is not
contained in $GL(n,\HH)$. Does this mean that quaternionic-K\"{a}hler
manifolds are not quaternionic? In reaction to this apparent contradiction,
S. Salamon defined quaternionic manifolds (see \cite{Sal}) as having a
holonomy reduction to $GL(n,\HH)\cdot Sp(1)$. Then
quaternionic-K\"{a}hler implies quaternionic, as you might expect. 

There are two interesting low dimensional coincidences in Berger's
list. The first $U(1) = SO(2)$ tells us the complex K\"{a}hler curves
are simply Riemannian surfaces. Moreover, because fixing a conformal structure
on a Riemannian surface corresponds to a holonomy reduction to
$\R^+SO(2)$, and $\R^+SO(2) = \C^* = GL(1,\C)$, geometrically
speaking, fixing a conformal structure is equivalent to fixing a
complex structure on two dimensional manifolds. 

The second coincidence $Sp(1)\cdot Sp(1) = SO(4)$ seems similar, with
``complex'' replaced by ``quaternionic''. We also have,
$$ GL(1,\HH)\cdot Sp(1) = \HH^*\cdot Sp(1) = \R^+Sp(1)\cdot Sp(1) =
\R^+SO(4) .$$
The implication is that fixing a quaternionic structure is equivalent
to fixing a conformal structure on four manifolds. But what exactly do
we mean by a ``quaternionic structure''?

\subsection{The Impasse}
The algebra of quaternions appears naturally as the generator of the
Brauer group of the reals, $\Br(\R) = \{\R,\HH\}$. The group
structure is given by the tensor product, moduli ``matrix''
algebra. It is not necessary to go into the details of the Brauer group
here, instead we simply note that $\HH$ generates $\Br(\R)$ because of
the following algebra isomorphism,
$$ \phi: \HH\otimes_\R\HH \simeq  \End_\R(\HH) ,$$
where $\phi(p\otimes q): v \mapsto p\cdot v\cdot q$ for any
$p,q,v\in\HH$. Note that we have used both the left and the right
module structures in defining $\phi$.

The Euclidean conformal group $\R^+SO(4)$ has a natural quaternionic
presentation using the isomorphism $\phi$. Let $i:
\HH\times\HH\to\HH\otimes\HH$ be the canonical map associated to the
tensor product. The image of the multiplicative group
$\HH^*\times\HH^*$ under these maps is precisely the conformal
group. We have the following exact sequence of groups,
$$\begin{CD} 1 @>>> \R^* @>>> \HH^*\times\HH^* @>{\phi\circ i}>>
  \R^+SO(4) @>>> 1 \end{CD}$$
where $\R^* \to \HH^*\times\HH^*$ acts as $r\mapsto (r,r^{-1})$. 

\begin{prop} The Euclidean conformal group in four dimensions appears
  in a natural and quaternionic way as,
$$ \R^+SO(4) = \{p\otimes q = i(p,q)\mid p,q \in\HH^* \}.$$ \end{prop}  

\textsc{Proof:} The Euclidean norm of $x\in\HH$ is $|x|^2 =
x\cdot\bar{x}$. Let $p\otimes q = i(p,q)$. Then,
\begin{align*} |p\otimes q (x)| &= |p\cdot x\cdot q| = \sqrt{p\cdot
    x\cdot q\cdot\bar{q}\cdot \bar{x}\cdot\bar{p}} \\ 
&= |p||q||x| = \lambda |x| \end{align*}
\Rbox

The above presentation of the conformal group, using the
isomorphism $\phi:\HH\otimes\HH\to\End(\HH)$, places equal emphasis
on the left and right module structures of $\HH$ on itself. Indeed,
the isomorphism $\phi$ is the $\HH$-bimodule structure on $\HH$! It is
then natural to consider the full bimodule structure as the important
structure that we want to integrate over four manifolds. However
this way is blocked.

\begin{prop} The automorphisms of $\HH$ considered as a $\HH$-bimodule
  are all scale multiples of the identity, 
$$ \Aut_{\HH\mbox{-bimodule}}(\HH) = \R^+\cdot \Id .$$ \end{prop}

\textsc{Proof:} This is simply a consequence of Shur's lemma applied
to the representations of $\M_4(\R)$.\Rbox

Thus a four manifold with an integrable $\HH$-bimodule structure
defined on the tangent bundle has holonomy contained in
$\R^+\cdot\Id$, which forces the manifold to be affine.

So we reach an impasse:
\begin{itemize} \item The Euclidean conformal group has a natural
  quaternionic presentation using the $\HH$-bimodule structure on
  $\HH$. 
\item The automorphisms of $\HH$ as a $\HH$-bimodule are simply scale
  multiples of the identity. \end{itemize}

We show that quaternionic gerbes provide a way of going past this
impasse. 

\subsection{The Suggested Solution}
The central idea is to use a more sophisticated way of ``gluing''
local data together.

Although $\HH$ has very few automorphisms when considered as an
$\HH$-bimodule, it does have an interesting group of automorphisms as
an $\R$-algebra,
$$ \Aut(\HH) = \Inn(\HH) = SO(3) .$$
Note that all the automorphisms are inner. The suggestion is to
consider the set of linear maps in $\End(\HH)$ that commute with the
$\HH$-bimodule structure, \textit{up to inner automorphisms}. Such a
map: $f: \HH\to \HH$ is required to satisfy the equation,
$$ f(p\cdot v\cdot q) = \alpha(p) f(v) \beta(q), $$
where $\alpha,\beta$ are inner automorphisms associated to $f$. It turns
out that the set of all such generalised automorphisms is precisely
the Euclidean conformal group, $i(\HH^*\times\HH^*)=\R^+SO(4)$.

The idea of allowing the two actions to commute \textit{up to an
automorphism} is natural in category theory. A gerbe is
a special kind of sheaf of categories and provides a rich enough
language to handle the inner automorphisms coherently.

An excellent presentation of the theory of \textit{Abelian} gerbes has
been presented by Jean-Luc Brylinski in ``Loop Spaces, Characteristic
Classes and Geometric Quantisation.'' \cite{Bry}. 

Nigel Hitchen, studying special Lagrangian sub-manifolds in dimension
three, has also made use of Abelian gerbes. Hitchen's approach
\cite{Hit} stresses the idea that Abelian gerbes certain cohomology
classes.  

Michael Murray has presented \cite{Mur} Abelian gerbes in a different
light as bundle gerbes. 

However the theory we are looking for is non-Abelian. L. Breen has
defined \textbf{non-Abelian} gerbes \cite{Br1,Br2} for arbitrary Lie
groups and has developed the theory of 2-gerbes. Breen's work
applies quite well to our present situation. 

\section{The Structure Groupoid}
Just as there is a structure group associated to a principal bundle, a
gerbe has an associated groupoid. In this section we will describe
the structure groupoid associated to a quaternionic gerbe.

Following Breen \cite{Br1}, we can associate to any \textit{crossed
module} a groupoid. We start from the crossed module defined by,
$$ \delta:\HH^*\to SO(3) .$$
where $\delta$ is the natural map onto the inner automorphisms,
$\delta(p) = p\otimes p^{-1}$, and $SO(3)$ acts on $\HH^*$ as
automorphisms.

Recently R. Brown and collaborators have been relating
groupoids and crossed modules to algebraic topology. (see \cite{Bro})

\begin{Def} The \textbf{quaternionic structure groupoid} $\cH$ is
  defined:
\begin{itemize} \item objects of $\cH$ are elements of $SO(3)$,
\item any element $p\in\HH^*$ is considered a morphism
  $p\in\cH(\alpha,\beta)$ 
$$ p : \alpha \to\beta $$
whenever $ \delta(p)\circ\alpha = \beta$. 
\item For any two morphisms $p:\alpha\to \beta$ and
  $q:\beta\to\gamma$, the composition is given by the map,
$$ q\circ p = qp : \alpha \to \gamma = \delta(qp)\alpha = \delta(q)\delta(p)\alpha.$$
\end{itemize}\end{Def}

It is easy to check that all of the axioms of a small category are
satisfied. In addition, because $\HH$ is a division algebra, all the
morphisms are invertible so that $\cH$ is a groupoid. 

Note that the set of all morphisms in $\cH$ consists of pairs
$(p,\alpha)\in\HH^*\times SO(3)$. We will abuse notation a little and
say that $\cH = \HH^*\times SO(3)$ as sets.

Although we have used the left $SO(3)$-action on itself, we have not
used the group structure on $SO(3)$. 

\subsection{Tensor product on $\cH$}
The small category $\cH$ has a monoidal structure coming
from the group structure on $SO(3)$,
$$\otimes : \cH \times \cH \to \cH .$$
We use the tensor product symbol because the central example of a
monoidal structure on a category is that of the tensor product in
vector spaces. This tensor product is not commutative, however it is
associative.  

For any $\alpha,\beta\in SO(3)$,
$$ \alpha\otimes\beta = \alpha\beta $$

For any two maps $p:\alpha\to\delta(p)\alpha$ and
$q:\beta\to\delta(q)\beta$ we define $p\otimes q$ as,
$$ p\otimes q = p\alpha[q] : \alpha\beta \to
\delta(p\alpha[q])\alpha\beta . $$
To see that $\otimes$ is well defined we should check that tensor
product of the ranges is the range of the tensor product of the maps,
\begin{align*} \delta(p\alpha[q])\alpha\beta &=
  \delta(p)\delta(\alpha[q])\alpha\beta =
  \delta(p)\alpha\delta(q)\alpha^{-1}\alpha\beta \\ 
&= (\delta(p)\alpha)\otimes(\delta(q)\beta) \end{align*}

Note that $\otimes$ is simply the semi-direct group structure coming from
the action of $SO(3)$ on $\HH^*$  via inner automorphisms,
$$ (\cH, \otimes) = \HH^* \rtimes SO(3) .$$
where $(p,\alpha)\otimes (q,\beta) = (p\alpha[q],\alpha\beta)$.

Moreover, this semi-direct product is isomorphic to the Euclidean
conformal group, 
$$ \HH^* \rtimes SO(3) \simeq \R^+SO(4) .$$

\subsection{$\HH$-bimodules} We can also represent the groupoid $\cH$ as a
category of quaternionic bimodules in such a way that the tensor
product is really a tensor product. In order to do this we need to
define carefully what we mean by an $\HH$-bimodule.

\begin{Def} An \textbf{$\HH$-bimodule} is a vector space $V$ with two
  commuting actions of the quaternions. Or equivalently, a bilinear map,
$$ \rho : \HH\times\HH \to \End(V) .$$ \end{Def}

Given an $\HH$-bimodule $(V,\rho)$ we can present the action as,
$$ \rho: \HH\otimes\HH \to \End(V),$$
by using the universal property of the tensor product.
In this form we see that the $\HH$-bimodules are simply the
modules over $\End(\HH)=\M_4(\R)$, the algebra of four by four
matrices. Therefore the only simple module is $\R^4$. For our purposes
we restrict ourselves to real four dimensional $\HH$-bimodules. 
The objects of $\cH$ will be identified with the four dimensional
$\HH$-bimodules.

Although all such objects are structurally identical, we will
distinguish between different quaternionic structures on the same
underlying vector space.

\begin{prop} Let $V$ be a $\HH$-bimodule in $\cH$. Then there is an
  $\HH$-bimodule isomorphism $\epsilon: \HH\to V$. \end{prop}

\textsc{Proof:} $V$ is a simple module over $\HH$, in two different
ways. Comparing these actions we define for each $v\in V$, $v\neq 0$,
a map $\HH\to\HH: p\mapsto p^v$ by the rule,
$$ p\cdot v = v\cdot p^v .$$
The map  $p\mapsto p^v$ is an $\R$-algebra automorphism for all $v$,
\begin{align*}  v \cdot (pq)^v &= (p q) \cdot v = p\cdot (q\cdot v) =
  p\cdot (v \cdot q^v) \\ 
&= (p\cdot v) \cdot q^v = v \cdot (p^vq^v) \end{align*}
So we have defined a map $V\to \Aut(\HH)=SO(3)$. To show
surjectivity we start by fixing some $v$ and define $\alpha[p] =
p^v$. For any $\beta$ in $\Aut(\HH)$ the tansitivity of $SO(3)$
implies that $\beta = \gamma\alpha$ for some $\gamma\in
SO(3)$. All automorphisms are inner, so there is some
$r\in\HH$ such that $\gamma = \delta(r) = r\otimes r^{-1}$. Then we
observe,
\begin{align*} p\cdot (v\cdot r^{-1}) &= (v \cdot p^v)\cdot r^{-1} =
  (v\cdot r^{-1}) \cdot (rp^vr^{-1}) \\ 
&= (v\cdot r^{-1}) \cdot \gamma[p^v] = (v\cdot r^{-1}) \cdot
\gamma\alpha[p] \\ 
&= (v\cdot r^{-1}) \cdot \beta[p] \end{align*}
and we see that $v\cdot r^{-1}$ maps onto $\beta$. It is easy to see
that the map $V\to SO(3)$ fibres through the projection $V\to P(V)$,
where $P(V)$ is the real projective space of one dimensional subspaces
in $V$. Let $e\in V$ be chosen in the preimage of the identity of
$SO(3)$. Then it is clear that $p\cdot e = e\cdot p$ for all
$p\in\HH$, and,
\begin{align*}\epsilon : \HH &\to V \\
 p &\to p\cdot e\end{align*} 
is the isomorphism of $\HH$-bimodules. \Rbox

We see from the above proof that each $\HH$-bimodule structure creates
an identification of $P(V)$ with $SO(3)$. Recall that as a smooth
manifold $SO(3)= \RP^3 = P(V)$. 

We can go in the other direction as well. Let $V$ be a right
$\HH$-module and $\alpha: V\to SO(3)$ be a $\HH^*$ equivariant map,
$$ \alpha(xp) = \delta(p)^{-1}\alpha(x).$$
Then we define a left $\HH^*$ action on $V$ as, 
$$ px = x\alpha(x)[p] $$
The left action commutes with the right action and so $V$ is an
$\HH$-bimodule,  
\begin{align*} p(xq) &= xq\alpha(xq)[p] = xq\delta(q)^{-1}\alpha(x)[p] \\
&= x\alpha(x)[p]q = (px)q \end{align*}

We distinguish the different objects in $\cH$ by using the different
identifications $P(V)\to SO(3)$ associated $\HH$-bimodule
structures. Two objects differ by an element of $SO(3)$. 

Now the tensor product can actually be represented as a tensor
product. If $V$ and $W$ are the $\HH$-bimodules associated to objects
$\alpha$ and $\beta$ in $\cH$, then the $\HH$-bimodule associated to
$\alpha\otimes\beta = \alpha\cdot \beta$ is $V\otimes_\HH W$.

A quaternionic gerbe consists of the structure groupoid fibred
over a four manifold. To see how we do that, we need a closer look at
the theory of sheaves of categories. 

\section{Sheaves of Categories or Stacks.}

A gerbe is a special kind of sheaf of categories. Our objective in
this section is to present enough of the general theory so that we can
understand what is the nature of gerbes, and how they can be
useful. We will not present a self contained account here, instead we
refer the reader to \cite{Bry}.

A presheaf of categories involves the interplay of locally defined
``objects'' and ``morphisms''. A \textbf{stack}\footnote{For us the
  terms ``stack'' and ``sheaf of categories'' refer to the same
  concept.} requires that the objects satisfy a 
descent property, \textit{up to an isomorphism}. The concept is
is quite flexible, but still very precise. The isomorphisms that glue
together the object data must satisfy additional coherence identities.

\subsection{Local Homeomorphisms}
Instead of working with the category of open sets on a
manifold $X$, we work with local homeomorphisms: continuous map $f :
Y \to X $ such that,
\begin{itemize}
\item any $y\in Y$ has an open neighbourhood $U$ whose image $f(U)$ is
open in $X$, and,
\item the restriction of $f$ to $U$ gives a homeomorphism between $U$
and $f(U)$.
\end{itemize}

\begin{Def} The category of \textbf{spaces over} $X$, $C_X$, has,
\begin{itemize}
\item objects are local homeomorphisms to $X$, $f:Y\to X,$
\item a morphism $g:(f:Y\to X)\to(h:Z\to X)$ is a local
  homeomorphisms, $ g : Y \to Z $ such that $f = h\circ g$. 
\end{itemize}\end{Def}

An important example to keep in mind is associated to an open cover
$\{U_i\}$ of $X$. The canonical projection $f:Y = \coprod_i U_i\to X$
from the disjoint union onto $X$ is a local homeomorphism. 

\subsection{Presheaves of Categories}
In the same way that a presheaf of sets is simply a functor from $C_X$
to the category of sets, a \textbf{presheaf of categories} over $X$ is a
functor $\cC$ from the category of spaces over $X$,
$C_X$, to the (bi)-category of small categories, functors and natural
transformations. Or, more explicitly,
\begin{itemize}\item to every local homeomorphism $f:Y\to X$ we
  associate a small category, 
$$ \cC(f:Y\to X)$$

\item to every arrow of local homeomorphisms $k:(Z,g)\to(Y,f)$ we
  associate a functor, 
$$ \cC(k) = k^{-1} : \cC(f:Y\to X) \to \cC(g:Z\to X),$$

\item to every composition $k\circ l:(W,h)\to
(Z,g)\to (Y,f)$ we associate an invertible natural transformation, 
$$ \theta_{k,l} : l^{-1}k^{-1} \Rightarrow (kl)^{-1} $$
\end{itemize}
This data must satisfy the following coherence condition,
\[ \begin{CD}
  m^{-1}l^{-1}k^{-1} @>{\theta_{k,l}}>> m^{-1}(lk)^{-1} \\
  @VV{\theta_{l,m}}V           @VV{\theta_{lk,m}}V \\
  (lm)^{-1}k^{-1} @>{\theta_{k,lm}}>> (lkm)^{-1} 
\end{CD}\]

It would be possible to define a presheaf of categories with the
requirement that $l^{-1}k^{-1}$ is strictly identical to
$(kl)^{-1}$. However that does not take advantage of the extra
flexibility provided. We will see later how quaternionic gerbes make
use of this flexibility.

\subsection{Descent for Morphisms}
Let $\cC$ be a presheaf of categories. We say that the
\textbf{morphisms satisfy descent} if for any two objects $A,B$ in
$\cC(f:Y\to X)$, the presheaf of sets on $Y$ defined by,
$$ \Hom(A,B)(k: Z\to Y) = \Hom(k^{-1}(A),k^{-1}(B)) $$
is actually a sheaf on $Y$.

We can explain this in terms of objects and maps more directly. Let
$V\subset X$ be a open neighbourhood, and let $A,B$ be objects in
$\cC(V) = \cC(V\hookrightarrow X)$. Now let $\{U_i\}$ be a open cover
of $V$. Take a collection of morphisms $\alpha_i: A|_{U_i} \to B|_{U_i}$,
where $\alpha_i\in\cC(U_i)$. The morphisms satisfy descent if,
$$ \alpha_i|_{U_{ij}} =  \alpha_j|_{U_{ij}} \quad \forall i,j,$$
implies the existence a unique morphism $\alpha:A\to B$ in
$\cC(V)$ such that $\alpha_i = \alpha|_{U_i}$  for all $i$.

In the above we have denoted $A|_{U_i}$ for the ``restriction'' of $A$
to $U_i$. Of course the restriction is really a functor
$\cC(V)\to\cC(V\cap U_i)$, and that functor is not necessarily trivial
or obvious. However the presentation becomes much easier to if we make
use of these small abuses of the notation.

\subsection{Descent for Objects}
The objects satisfy a much more complicated descent property, making
use of the natural transformations appearing in the definition.

Let $\cC$ be a presheaf of categories. Let $V$ be any open
set in $X$ and $f:Y\to V$ be any surjective local homeomorphism. 
The \textbf{descent} data for any $A\in\cC(Y)$ consists of an isomorphism
$\phi:p_2^{-1}(A) \to p_1^{-1}(A)$ in $\cC(Y\times_XY)$ such that,
$$ p_{12}^{-1}(\phi)\circ p_{23}^{-1}(\phi)\circ p_{31}^{-1}(\phi) =
\Id_{p^{-1}_1(A)} $$ 
in $\cH(Y\times_XY\times_XY)$.\footnote{The natural
  projections $Y\times_XY \to Y$ are denoted  $p_1$ and $p_2$, and the
  three natural projections $Y\times_XY\times_XY\to Y\times_XY$ are denoted by
  $p_{12}$,$p_{23}$ and $p_{13}$.}

We say that the \textbf{objects satisfy descent} if every pair
$(A,\phi)$ as above implies the existence of an object $A'\in \cC(V)$
and an isomorphism $\psi:f^{-1}(A') \to A$ in $\cC(Y)$ such that the
following diagram in $\cC(Y\times_XY)$ commutes, 
\[\begin{CD}
 p_1^{-1}f^{-1}(A') @>{\theta^{-1}_{f,p_2}\theta_{f,p_1}}>>
 p_2^{-1}f^{-1}(A') \\
 @V{\psi}VV   @V{\psi}VV \\
 p_1^{-1}(A) @>{\phi}>> p_2^{-1}(A) 
\end{CD}\]

This rather complicated prescription can also be understood in terms of
open sets in the normal sense. 

Let $\{U_i\}$ be an open cover of $V\subset X$. The descent data is
equivalent to a set of local objects $A_i\in\cC(U_i)$, and isomorphisms
$\phi_{ij} : A_i|_{U_{ij}} \to A_j|_{U_{ij}}$ in $\cC(U_i\cap
U_j)$. The isomorphisms are required to satisfy,
$$ \phi_{ik}|_{U_{ijk}} = \phi_{ij}|_{U_{ijk}}\circ\phi_{jk}|_{U_{ijk}}, $$
in the category over the triple intersection, $\cC(U_i\cap U_j\cap
U_k)$. Again note that we have implicitly used the natural
transformations by glossing over the restrictions. 

\begin{Def} A \textbf{stack} (or sheaf of categories) on $X$ is a
  presheaf of categories where objects and morphisms satisfy the
  descent conditions above. \end{Def}

\section{Quaternionic Gerbes}
Now we refine the notion of a stack to that of a gerbe by imposing
three more conditions:
\begin{enumerate}
\item gerbes take values in groupoids, the full
  sub-category of small categories whose morphisms are invertible. 
\item gerbes are \textbf{locally non-empty}. This means that there
  exists a surjective local homeomorphism $f:Y\to X$ such that
  $\cC(Y)$ is non-empty. We could also state this by saying that there
  exists an covering $\{U_i\}$ of $X$ such that the $\cC(U_i)$ are all
  non-empty. 
\item gerbes are \textbf{locally connected}. This means that for any
  two objects $A,B$ in $\cC(f:Y\to X)$, there exists an surjective
  local homeomorphism $g:Z\to Y$ such that $g^{-1}(A)$ and $g^{-1}(B)$
  are isomorphic. In terms of covers: if $A,B$ are objects in $\cC(U)$
  for some $U\subset X$, then there exists an open covering $\{U_i\}$
  of $U$ such that $A\mid_{U_i}$ is isomorphic to $B\mid_{U_i}$ for
  all $i$.
\end{enumerate}

\begin{Def} A \textbf{gerbe} on $X$ is a locally non-empty and locally
connected sheaf of groupoids on $X$. \end{Def}

For any group $G$ let $\underline{G}_X$ be the sheaf of $G$-valued
functions on $X$. A gerbe is said to have \textbf{band}
in $G$ if for any object $A\in \cC(f:Y\to X)$, the sheaf
$\Aut(A)$ of automorphisms of $A$ on $Y$ is isomorphic to
$\underline{G}_Y$, and the isomorphism $\alpha:\Aut(A) \to
\underline{G}_Y$ is unique up to an inner automorphism of $G$.

\begin{Def} \textbf{Quaternionic Gerbe} is a gerbe with band in
  $\HH^*$. \end{Def}

\subsection{Neutral Gerbes.}
A gerbe $\cG$ is said to be \textbf{neutral} if there exists a global
object, $A\in\cG(X)$. Because the automorphism sheaf of $\Aut(A)$ is
isomorphic to the sheaf of $\HH^*$-valued functions, we
can identify the groupoid $\cG(X)$ with the groupoid of principal
$\HH^*$-bundles,
\begin{align*} \phi : \cG(X) &\to \Tor(\Aut(A)) \\
B &\mapsto \Isom(B\to A) \end{align*}

Quaternionic gerbes are locally non-empty so we can always find
an object $A_U\in\cG(U)$ over the open set $U$. Using that local
object we can identify $\cG(U)$ with the groupoid of $\HH^*$-bundles over
$U$. Local non-emptiness implies that the gerbe is locally neutral.

In order to understand this local neutrality, it is helpful to
consider an analogy with the relation between a principal $G$-bundle
and an associated vector bundle. To any vector bundle we can associate
the principal bundle of frames. The local neutralisation associated to
a local object $A_U$ is sort of ``frame'' for $\cG$ over $U$. The set
of all frames for $\cG$ forms a local groupoid. 

Assuming that $U$ is contractable, all the $\HH^*$-bundles differ by
an automorphisms valued function $\alpha:U\to SO(3)$. We define a
$\HH^*$-bundle associated to $A_U$ and $\alpha$ by letting $A_U^\alpha
= A_U$ as a fibre bundle. The action of $\HH^*$ however is twisted by
$\alpha$. Let $a\in A_U$ and let $a^\alpha$ be the same element
considered in $A_U^\alpha$. Then for any quaternion $p\in\HH^*$,
$$ a^\alpha \cdot p = (a \cdot \alpha[p])^\alpha .$$

If $U$ is not contractable there can be topologically inequivalent
$\HH^*$-bundles. Then we can replace the function $\alpha$ above with
an $\HH$-bimodule $M\to U$. If $A_U$ and $B_U$ are two different
objects in $\cG(U)$, then there is an $\HH$-bimodule $M$ such that,
$$ A\otimes_\HH M = B .$$

The local groupoid $\cH(U)$ consists of the ``frames'' of
$\cG(U)$. Note that because of the local neutrality axiom, all
quaternionic gerbes look the same locally.

\subsection{The Local Groupoid}
We describe here the local structure of $\cH$.

An \textbf{object} of the local groupoid $\cH(U)$ is a diagram of
the form,
\[\begin{CD} A @>{\alpha}>> \Aut(\HH) \\
 @V{\pi}VV \\
 U \end{CD} \]
where $\pi:A\to U$ is a principle $\HH^*$-bundle and $\alpha$ is an
$\HH^*$-equivariant map $\alpha(xp) = \delta(p)^{-1}\alpha(x)$. 

As we have seen, this data can also be presented in terms of
$\HH$-bitorsors. 

An $\HH^*$-bitorsor is a principle right $\HH^*$-bundle
that is also a principle left $\HH^*$-bundle for a commuting action
of $\HH^*$. For any  $(A,\alpha) \in \cH(U)$, the left $\HH^*$ action
on $A$ is, $ px = x\alpha(x)[p]$.

The morphisms of $\HH$-bitorsors are simply bundle maps that commute
with both the left and right actions.

\subsection{Tensor Product on $\cH(U)$}
In terms of bitorsors we can present the product structure on 
$\cH(U)$ by using the quaternionic tensor product. 

For any $A, B\in\cH(U)$,
$$ A\otimes_\HH B = A\otimes_\R B / \sim $$
where $x p\otimes y \sim x\otimes py$.

Assuming $U$ is contractable and by fixing a coordinate basis, we
get a canonical trivialisation of the tangent bundle, $TU =
U\times\HH$. In this way $TU$ can be considered as an object in
$\cH(U)$.

Relative to this fixed object, all the others are given by $SO(3)$
valued functions on $U$, the morphisms are given by $\HH^*$ valued
functions. 

Over $U\subset\HH$ the local groupoid consists of sections
$C^\infty(U,\cH)$. However the strength of this approach is in terms
of the global structure. A global quaternionic gerbe is given in terms
of ``transition functions''.

\subsection{Transition Functions or Bitorsor Cocyle.}
The transition functions for a quaternionic gerbe are given in
terms of $\HH$-bitorsors. Maybe we should say ``transition bitorsors''.

Let $\cG$ be a quaternionic gerbe on $X$ and let $\{U_i\}$ be a good
cover.\footnote{All intersections $U_i\cap U_j$ are contractable.} Choose
$A_i\in\cG(U_i)$. Then $\Aut(A_i)$ is isomorphic to the sheaf of
$\HH^*$-valued functions. Using $A_i$ we have the following local
neutralisation,   
\begin{align*} \Phi_i : \cG(U_i) &\to \Tor(\Aut(A_i)) \\
B &\mapsto \Isom(B\to A_i) \end{align*}

On any intersection $U_{ij}=U_i\cap U_j$ we can define,
$$ E_{ij} = \Isom(A_j\mid_{U_{ij}}, A_i\mid_{U_{ij}}),$$

The $E_{ij}$ are $\HH$-bitorsors and are the transition
functions. The two $\HH$-actions are given by the composition of an
isomorphism with automorphisms of $A_i|_U{_{ij}}$ and $A_j|_U{_{ij}}$,
which are each isomorphic to $\HH^*$-valued functions. Note that the
isomorphisms $\underline{\HH}\simeq \Aut(A_i)$ are unique up to an
automorphism. To be really careful we should take care of those
automorphisms as well, however that will work will be presented in a
comprehensive way later.

The $\HH$-bitorsors $E_{ij}$ need to be compared over triple
intersections. The natural transformations in the definition of a
stack give us the following morphisms as extra data:
$$ \psi_{ijk} : E_{ij}\otimes_\HH E_{jk} \to E_{ik},$$

These morphisms live in $\cH(U_{ijk})$ and must satisfy the following
coherence condition on four intersections,
\[\begin{CD}
E_{ij}\otimes_\HH E_{jk}\otimes_\HH E_{kl} @>{\psi_{ijk}\otimes\Id}>> 
E_{ik}\otimes_\HH E_{kl} \\
@V{\Id\otimes\psi_{jkl}}VV  @VV{\psi_{ikl}}V \\
E_{ij}\otimes_\HH E_{jl} @>{\psi_{ijl}}>> E_{il} 
\end{CD} \]

The pair $(E_{ij},\psi_{ijk})$ is called a \textbf{quaternionic
  bitorsor cocycle} on $X$.

Of course the above description of a particular quaternionic gerbe
depends on the choice of $A_i\in \cG(U_i)$. We can measure the dependence
on those choices with a coboundary.

\subsection{Coboundary}

Let $B_i$ be a different choice of local objects and
$(F_{ij},\phi_{ijk})$ be the associated bitorsor cocyle.

Let $M_i\in\cH(U_i)$ be defined by,
$$ B_i = A_i\otimes_\HH M_i $$

The pair $(M_i, \nu_{ij})$ is a \textbf{coboundary}
relating $(F_{ij},\phi_{ijk})$ to $(E_{ij},\psi_{ijk})$ if $\nu_{ij}$
is a map in $\cH(U_{ij})$,
$$ \nu_{ij} : F_{ij} \to M^\circ_i\otimes_\HH E_{ij}\otimes_\HH M_j $$
such that as morphisms in $\cH(U_{ijk})$,
$$ \nu_{ik}\circ\phi_{ijk} = \psi_{ijk}\circ(\nu_{ij}\otimes\nu_{jk}) $$

We can present this equation with a commutative diagram,
\[\begin{CD}
F_{ij}\otimes_\HH F_{jk} @>{\nu_{ij}\otimes\nu_{jk}}>>
M^\circ_i\otimes_\HH E_{ij}\otimes_\HH E_{jk}\otimes_\HH M_k \\
@V{\phi_{ijk}}VV   @VV{\Id\cdot\psi_{ijk}\cdot\Id}V \\
F_{ik} @>>{\nu_{ik}}> M^\circ_i\otimes_\HH E_{ik}\otimes_\HH M_k
\end{CD}\]

It was demonstrated in \cite{Fin} that coboundaries define an
equivalence relation on the set of 
quaternionic bitorsor cocycles. Moreover, it is possible to construct
a quaternionic gerbe from a given cocyle, and that gerbe will be
isomorphic to any gerbe constructed from a cocyle from the same
equivalence class.

Although we have used the terminology of cohomology at present there
is no actual theory of $\HH$-valued cohomology. We use the terminology
because it is convenient, and perhaps to be a little optimistic.

\section{Conformal Four Manifolds}

A conformal structure on a four manifold is a reduction of the frame
bundle to $\R^+SO(4)$. As we saw at the beginning, the Euclidean
conformal group can be presented using the groupoid $\cH$ with its
tensor product acting as the group structure.

We will indicate here briefly how to construct a quaternionic bitorsor
cocycle from a given a conformal structure on a four manifold $X$. The
presentation here is very sketchy and a more detailed presentation is
being prepared.

We can choose charts $\{\psi_i : U_i \to \HH\}$
that are compatible with the conformal structure, i.e.,
$$ \partial(\psi_i\circ\psi_j^{-1}) \in i(\HH^*\times\HH^*),$$
where $\partial(f)$ is the Jacobian matrix of $f$ considered as an
element of $\HH\otimes\HH$. Therefore there are $\HH^*$-valued
functions $x_{ij}$ and $y_{ij}$ on $U_{ij}$ such that,
$$ \partial(\psi_i\circ\psi_j^{-1}) = x_{ij}\otimes y_{ji} .$$

Over each chart $U_i$ the tangent bundle has a canonical $\HH$-bitorsor
structure coming from the coordinate $\psi$. The \textbf{tangent
  gerbe cocycle} allows us to relate these various $\HH$-bitorsor
structures. 

The $x_{ij}\otimes y_{ij}$ can be used to define $E_{ij}$ by twisting
the left and right $\HH$-actions by $\delta(x_{ij})$ and
$\delta(y_{ij})$ respectively. In terms of an $SO(3)$ valued function,
we can define $E_{ij}$ relative to $TU_i$ with the function
$\delta(y_{ij}x_{ij})$. 

Over the triple intersections $U_{ijk}$ it is possible to construct
isomorphisms $\phi_{ijk}: E_{ij}\otimes E_{jk} \to E_{ik}$.

It can also be shown that the $\psi_i$ are coordinate charts compatible
with the conformal structure if and only if the $(E_{ij},\phi_{ijk})$
form a quaternionic gerbe cocyle.

\end{document}